\documentclass[12pt]{article}
\usepackage{theorem}
\usepackage{amssymb}
\usepackage{graphicx}
\usepackage[mathscr]{eucal}
\usepackage{amsbsy}
\usepackage{amsmath}
\usepackage{cite}
\newtheorem{prop}{}[section]

{\theorembodyfont{\upshape} \newtheorem{rema}[prop]{}}
\newcommand{\boma}[1]{{\mbox{\boldmath $#1$} }}

\begin{document}
\newcommand{\uper}[1]{\stackrel{\barray{c} {~} \\ \mbox{\footnotesize{#1}}\farray}{\longrightarrow} }
\newcommand{\nop}[1]{ \|#1\|_{\piu} }
\newcommand{\no}[1]{ \|#1\| }
\newcommand{\nom}[1]{ \|#1\|_{\meno} }
\newcommand{\UU}[1]{e^{#1 \AA}}
\newcommand{\UD}[1]{e^{#1 \Delta}}
\newcommand{\bb}[1]{\mathbb{{#1}}}
\newcommand{\HO}[1]{\bb{H}^{{#1}}}
\newcommand{\Hz}[1]{\bb{H}^{{#1}}_{\zz}}
\newcommand{\Hs}[1]{\bb{H}^{{#1}}_{\ss}}
\newcommand{\Hg}[1]{\bb{H}^{{#1}}_{\gg}}
\newcommand{\HMs}[1]{\bb{H}^{{#1}}}
\newcommand{\CMs}[1]{\bb{C}^{{#1}}}
\newcommand{\EMs}[1]{\bb{E}^{{#1}}}
\newcommand{\EMnus}[1]{\bb{E}^{{#1}}_{\nu}}
\newcommand{\HM}[1]{\bb{H}^{{#1}}_{\so}}
\newcommand{\CM}[1]{\bb{C}^{{#1}}_{\so}}
\newcommand{\EM}[1]{\bb{E}^{{#1}}_{\so}}
\newcommand{\EMnu}[1]{\bb{E}^{{#1}}_{\so\,\nu}}
\def\uz{u_{0}}
\def\vz{v_{0}}
\def\vi{v}
\def\ef{\psi}
\def\fun{\mathcal{F}}
\def\fun{{\tt f}}
\def\tvainf{\vspace{-0.4cm} \barray{ccc} \vspace{-0,1cm}{~}
\\ \vspace{-0.2cm} \longrightarrow \\ \vspace{-0.2cm} \scriptstyle{T \vain + \infty} \farray}
\def\De{F}
\def\er{\epsilon}
\def\erd{\er_0}
\def\Tn{T_{\star}}
\def\Tc{T_{\tt{c}}}
\def\Tb{T_{\tt{b}}}
\def\Tl{\mathscr{T}}
\def\Tm{T}
\def\pa{p_{\tt{a}}}
\def\Ta{T_{\tt{a}}}
\def\ua{u_{\tt{a}}}
\def\Tg{T_{G}}
\def\Tgg{T_{I}}
\def\Tw{T_{\star}}
\def\Ts{T_{\Ss}}
\def\Tr{\Tl}
\def\Sp{\Ss'}
\def\Tsp{T_{\Sp}}
\def\vsm{\vspace{-0.1cm}\noindent}
\def\comple{\scriptscriptstyle{\complessi}}
\def\ug{u_G}
\def\nume{0.407}
\def\numerob{0.00724}
\def\deln{7/10}
\def\delnn{\dd{7 \over 10}}
\def\e{c}
\def\p{p}
\def\z{z}
\def\symd{{\mathfrak S}_d}
\def\del{\omega}
\def\Del{\delta}
\def\Di{\Delta}
\def\Ss{{\mathscr{S}}}
\def\Ww{{\mathscr{W}}}
\def\mmu{\hat{\mu}}
\def\rot{\mbox{rot}\,}
\def\curl{\mbox{curl}\,}
\def\Mm{\mathscr M}
\def\XS{\boma{x}}
\def\TS{\boma{t}}
\def\Lam{\boma{\eta}}
\def\DS{\boma{\rho}}
\def\KS{\boma{k}}
\def\LS{\boma{\lambda}}
\def\PR{\boma{p}}
\def\VS{\boma{v}}
\def\ski{\! \! \! \! \! \! \! \! \! \! \! \! \! \!}
\def\h{L}
\def\EMP{M'}
\def\R{R}
\def\Aa{{\mathscr{A}}}
\def\Rr{{\mathscr{R}}}
\def\Zz{{\mathscr{Z}}}
\def\Jj{{\mathscr{J}}}
\def\E{E}
\def\FFf{\mathscr{F}}
\def\A{F}
\def\Xim{\Xi_{\meno}}
\def\Ximn{\Xi_{n-1}}
\def\lan{\lambda}
\def\om{\omega}
\def\Om{\Omega}
\def\Sim{\Sigm}
\def\Sip{\Delta \Sigm}
\def\Sigm{{\mathscr{S}}}
\def\Ki{{\mathscr{K}}}
\def\Hi{{\mathscr{H}}}
\def\zz{{\scriptscriptstyle{0}}}
\def\ss{{\scriptscriptstyle{\Sigma}}}
\def\gg{{\scriptscriptstyle{\Gamma}}}
\def\so{\ss \zz}
\def\Dv{\bb{\DD}'}
\def\Dz{\bb{\DD}'_{\zz}}
\def\Ds{\bb{\DD}'_{\ss}}
\def\Dsz{\bb{\DD}'_{\so}}
\def\Dg{\bb{\DD}'_{\gg}}
\def\Ls{\bb{L}^2_{\ss}}
\def\Lg{\bb{L}^2_{\gg}}
\def\bF{{\bb{V}}}
\def\Fz{\bF_{\zz}}
\def\Fs{\bF_\ss}
\def\Fg{\bF_\gg}
\def\Pre{P}
\def\UUU{{\mathcal U}}
\def\fiapp{\phi}
\def\PU{P1}
\def\PD{P2}
\def\PT{P3}
\def\PQ{P4}
\def\PC{P5}
\def\PS{P6}
\def\Q{P6}
\def\X{Q2}
\def\Xp{Q3}
\def\Vi{V}
\def\bVi{\bb{V}}
\def\K{V}
\def\Ks{\bb{\K}_\ss}
\def\Kz{\bb{\K}_0}
\def\KM{\bb{\K}_{\, \so}}
\def\HGG{\bb{H}^\G}
\def\HG{\bb{H}^\G_{\so}}
\def\EG{{\mathfrak{P}}^{\G}}
\def\G{G}
\def\de{\delta}
\def\esp{\sigma}
\def\dd{\displaystyle}
\def\LP{\mathfrak{L}}
\def\dive{\mbox{div}}
\def\la{\langle}
\def\ra{\rangle}
\def\um{u_{\meno}}
\def\uv{\mu_{\meno}}
\def\Fp{ {\textbf F_{\piu}} }
\def\Ff{ {\textbf F} }
\def\Fm{ {\textbf F_{\meno}} }
\def\Eb{ {\textbf E} }
\def\piu{\scriptscriptstyle{+}}
\def\meno{\scriptscriptstyle{-}}
\def\omeno{\scriptscriptstyle{\ominus}}
\def\Tt{ {\mathscr T} }
\def\Xx{ {\textbf X} }
\def\Yy{ {\textbf Y} }
\def\Ee{ {\textbf E} }
\def\VP{{\mbox{\tt VP}}}
\def\CP{{\mbox{\tt CP}}}
\def\cp{$\CP(f_0, t_0)\,$}
\def\cop{$\CP(f_0)\,$}
\def\copn{$\CP_n(f_0)\,$}
\def\vp{$\VP(f_0, t_0)\,$}
\def\vop{$\VP(f_0)\,$}
\def\vopn{$\VP_n(f_0)\,$}
\def\vopdue{$\VP_2(f_0)\,$}
\def\leqs{\leqslant}
\def\geqs{\geqslant}
\def\mat{{\frak g}}
\def\tG{t_{\scriptscriptstyle{G}}}
\def\tN{t_{\scriptscriptstyle{N}}}
\def\TK{t_{\scriptscriptstyle{K}}}
\def\CK{C_{\scriptscriptstyle{K}}}
\def\CN{C_{\scriptscriptstyle{N}}}
\def\CG{C_{\scriptscriptstyle{G}}}
\def\CCG{{\mathscr{C}}_{\scriptscriptstyle{G}}}
\def\tf{{\tt f}}
\def\ti{{\tt t}}
\def\ta{{\tt a}}
\def\tc{{\tt c}}
\def\tF{{\tt R}}
\def\C{{\mathscr C}}
\def\P{{\mathscr P}}
\def\V{{\mathscr V}}
\def\TI{\tilde{I}}
\def\TJ{\tilde{J}}
\def\Lin{\mbox{Lin}}
\def\Hinfc{ H^{\infty}(\reali^d, \complessi) }
\def\Hnc{ H^{n}(\reali^d, \complessi) }
\def\Hmc{ H^{m}(\reali^d, \complessi) }
\def\Hac{ H^{a}(\reali^d, \complessi) }
\def\Dc{\DD(\reali^d, \complessi)}
\def\Dpc{\DD'(\reali^d, \complessi)}
\def\Sc{\SS(\reali^d, \complessi)}
\def\Spc{\SS'(\reali^d, \complessi)}
\def\Ldc{L^{2}(\reali^d, \complessi)}
\def\Lpc{L^{p}(\reali^d, \complessi)}
\def\Lqc{L^{q}(\reali^d, \complessi)}
\def\Lrc{L^{r}(\reali^d, \complessi)}
\def\Hinfr{ H^{\infty}(\reali^d, \reali) }
\def\Hnr{ H^{n}(\reali^d, \reali) }
\def\Hmr{ H^{m}(\reali^d, \reali) }
\def\Har{ H^{a}(\reali^d, \reali) }
\def\Dr{\DD(\reali^d, \reali)}
\def\Dpr{\DD'(\reali^d, \reali)}
\def\Sr{\SS(\reali^d, \reali)}
\def\Spr{\SS'(\reali^d, \reali)}
\def\Ldr{L^{2}(\reali^d, \reali)}
\def\Hinfk{ H^{\infty}(\reali^d, \KKK) }
\def\Hnk{ H^{n}(\reali^d, \KKK) }
\def\Hmk{ H^{m}(\reali^d, \KKK) }
\def\Hak{ H^{a}(\reali^d, \KKK) }
\def\Dk{\DD(\reali^d, \KKK)}
\def\Dpk{\DD'(\reali^d, \KKK)}
\def\Sk{\SS(\reali^d, \KKK)}
\def\Spk{\SS'(\reali^d, \KKK)}
\def\Ldk{L^{2}(\reali^d, \KKK)}
\def\Knb{K^{best}_n}
\def\sc{\cdot}
\def\k{\mbox{{\tt k}}}
\def\x{\mbox{{\tt x}}}
\def\g{ {\textbf g} }
\def\QQQ{ {\textbf Q} }
\def\AAA{ {\textbf A} }
\def\gr{\mbox{gr}}
\def\sgr{\mbox{sgr}}
\def\loc{\mbox{loc}}
\def\PZ{{\Lambda}}
\def\PZAL{\mbox{P}^{0}_\alpha}
\def\epsilona{\epsilon^{\scriptscriptstyle{<}}}
\def\epsilonb{\epsilon^{\scriptscriptstyle{>}}}
\def\lgraffa{ \mbox{\Large $\{$ } \hskip -0.2cm}
\def\rgraffa{ \mbox{\Large $\}$ } }
\def\restriction{\upharpoonright}
\def\M{{\scriptscriptstyle{M}}}
\def\m{m}
\def\Fre{Fr\'echet~}
\def\I{{\mathcal N}}
\def\ap{{\scriptscriptstyle{ap}}}
\def\fiap{\varphi_{\ap}}
\def\dfiap{{\dot \varphi}_{\ap}}
\def\DDD{ {\mathfrak D} }
\def\BBB{ {\textbf B} }
\def\EEE{ {\textbf E} }
\def\GGG{ {\textbf G} }
\def\TTT{ {\textbf T} }
\def\KKK{ {\textbf K} }
\def\HHH{ {\textbf K} }
\def\FFi{ {\bf \Phi} }
\def\GGam{ {\bf \Gamma} }
\def\sc{ {\scriptstyle{\bullet} }}
\def\a{a}
\def\ep{\epsilon}
\def\c{\kappa}
\def\parn{\par \noindent}
\def\teta{M}
\def\elle{L}
\def\ro{\rho}
\def\al{\alpha}
\def\si{\sigma}
\def\be{\beta}
\def\ga{\gamma}
\def\te{\vartheta}
\def\ch{\chi}
\def\et{\eta}
\def\complessi{{\bf C}}
\def\len{{\bf L}}
\def\reali{{\bf R}}
\def\interi{{\bf Z}}
\def\Z{{\bf Z}}
\def\naturali{{\bf N}}
\def\Sfe{ {\bf S} }
\def\To{ {\bf T} }
\def\Td{ {\To}^d }
\def\Tt{ {\To}^3 }
\def\Zd{ \interi^d }
\def\Zt{ \interi^3 }
\def\Zet{{\mathscr{Z}}}
\def\Ze{\Zet^d}
\def\T1{{\textbf To}^{1}}
\def\es{s}
\def\ee{{E}}
\def\FF{\mathcal F}
\def\FFu{ {\textbf F_{1}} }
\def\FFd{ {\textbf F_{2}} }
\def\GG{{\mathcal G} }
\def\EE{{\mathcal E}}
\def\KK{{\mathcal K}}
\def\PP{{\mathcal P}}
\def\PPP{{\mathscr P}}
\def\PN{{\mathcal P}}
\def\PPN{{\mathscr P}}
\def\QQ{{\mathcal Q}}
\def\J{J}
\def\Np{{\hat{N}}}
\def\Lp{{\hat{L}}}
\def\Jp{{\hat{J}}}
\def\Pp{{\hat{P}}}
\def\Pip{{\hat{\Pi}}}
\def\Vp{{\hat{V}}}
\def\Ep{{\hat{E}}}
\def\Gp{{\hat{G}}}
\def\Kp{{\hat{K}}}
\def\Ip{{\hat{I}}}
\def\Tp{{\hat{T}}}
\def\Mp{{\hat{M}}}
\def\La{\Lambda}
\def\Ga{\Gamma}
\def\Si{\Sigma}
\def\Upsi{\Upsilon}
\def\Gam{\Gamma}
\def\Gag{{\check{\Gamma}}}
\def\Lap{{\hat{\Lambda}}}
\def\Upsig{{\check{\Upsilon}}}
\def\Kg{{\check{K}}}
\def\ellp{{\hat{\ell}}}
\def\j{j}
\def\jp{{\hat{j}}}
\def\BB{{\mathcal B}}
\def\LL{{\mathcal L}}
\def\MM{{\mathcal U}}
\def\SS{{\mathcal S}}
\def\DD{D}
\def\Dd{{\mathcal D}}
\def\VV{{\mathcal V}}
\def\WW{{\mathcal W}}
\def\OO{{\mathcal O}}
\def\RR{{\mathcal R}}
\def\TT{{\mathcal T}}
\def\AA{{\mathcal A}}
\def\CC{{\mathcal C}}
\def\JJ{{\mathcal J}}
\def\NN{{\mathcal N}}
\def\HH{{\mathcal H}}
\def\XX{{\mathcal X}}
\def\XXX{{\mathscr X}}
\def\YY{{\mathcal Y}}
\def\ZZ{{\mathcal Z}}
\def\CC{{\mathcal C}}
\def\cir{{\scriptscriptstyle \circ}}
\def\circa{\thickapprox}
\def\vain{\rightarrow}
\def\salto{\vskip 0.2truecm \noindent}
\def\spazio{\vskip 0.5truecm \noindent}
\def\vs1{\vskip 1cm \noindent}
\def\fine{\hfill $\square$ \vskip 0.2cm \noindent}
\def\ffine{\hfill $\lozenge$ \vskip 0.2cm \noindent}
\newcommand{\rref}[1]{(\ref{#1})}
\def\beq{\begin{equation}}
\def\feq{\end{equation}}
\def\beqq{\begin{eqnarray}}
\def\feqq{\end{eqnarray}}
\def\barray{\begin{array}}
\def\farray{\end{array}}
\makeatletter \@addtoreset{equation}{section}
\renewcommand{\theequation}{\thesection.\arabic{equation}}
\makeatother
\begin{titlepage}
{~}
\vspace{1cm}
\begin{center}
{\Large \textbf{Smooth solutions of the Euler and }}
\vskip 0.2cm
{\Large \textbf{Navier-Stokes equations from the a posteriori}}
\vskip 0.2cm
{\Large \textbf{analysis of approximate solutions}}
\end{center}
\vspace{0.5truecm}
\begin{center}
{\large
Carlo Morosi$\,{}^a$, Livio Pizzocchero$\,{}^b$({\footnote{Corresponding author}})} \\
\vspace{0.5truecm} ${}^a$ Dipartimento di Matematica, Politecnico di Milano,
\\ P.za L. da Vinci 32, I-20133 Milano, Italy \\
e--mail: carlo.morosi@polimi.it \\
${}^b$ Dipartimento di Matematica, Universit\`a di Milano\\
Via C. Saldini 50, I-20133 Milano, Italy\\
and Istituto Nazionale di Fisica Nucleare, Sezione di Milano, Italy \\
e--mail: livio.pizzocchero@unimi.it
\end{center}
\begin{abstract}
The main result of (C. Morosi and L. Pizzocchero, Nonlinear Analysis, 2012 \cite{appeul}) is
presented in a variant, based on a $C^\infty$ formulation of the
Cauchy problem; in this approach, the a posteriori analysis of an approximate solution
gives a bound on the Sobolev distance of any order between the exact and the
approximate solution.
\end{abstract}
\vspace{1cm} \noindent
\hfill \par
\par \vspace{0.05truecm} \noindent \textbf{AMS 2000 Subject classifications:} 35Q30, 76D03, 76D05.
\end{titlepage}
\section{Introduction}
Let us consider the homogeneous incompressible Navier-Stokes (NS) equations on a torus
$\Td$ of arbitrary dimension $d$; these read
\beq {\partial u \over \partial t}  = \nu \Delta u + \PPP(u,u) + f.
\label{eulnu} \feq
Here: $u = u(x,t)$ is the divergence free velocity field, depending on $x \in \Td$ and
on time $t$; $\nu \geqs 0$ is the kinematic viscosity, and $\Delta$ is the Laplacian
of $\Td$; $\PPP$ is the bilinear map sending any two sufficiently regular
vector fields $v, w : \Td \vain \reali^d$ into
\beq \PPP(v,w) := - \LP(v \sc \partial w)~. \label{defp} \feq
In the above $v \sc \partial w: \Td \vain \reali^d$ is
the vector field of components $(v \sc \partial w)_r = \sum_{s=1}^d v_s \partial_s w_r$,
and $\LP$ is the Leray projection onto the space of divergence free vector fields. Finally,
Eq.\,\rref{eulnu} contains the (Leray projected) density $f = f(x,t)$ of
the external forces. In the inviscid case $\nu=0$, the NS equations become the Euler equations.
\par
In our previous work \cite{appeul}, exact and approximate
solutions of the NS Cauchy problem have been discussed in a
framework based on the Sobolev spaces $\HM{n}$, for appropriate
$n$. For each real $n$, $\HM{n}$ consists of the (distributional)
vector fields $v: \Td \vain \reali^d$ with vanishing divergence
and mean such that $(-\Delta)^{n/2} v $ is square integrable;
this space carries the inner product
$\la v | w \ra_n := \la (-\Delta)^{n/2} v |  (-\Delta)^{n/2} w \ra_{L^2}$
and the corresponding norm $\|~\|_n$
(see the forthcoming
Eqs.\,\rref{defhn}\rref{definner}). After fixing an order $n > d/2
+ 1$, a forcing $f$ and an initial datum $u_0 \in \HM{n+2}$, in
the cited work we have discussed exact and approximate solutions
of the NS Cauchy problem in the functional class $C([0,T),
\HM{n+2}) \cap C^1([0,T), \HM{n})$ (with $T$ possibly depending on
the solution; in the Euler case $\nu=0$, one can harmlessly
replace $\HM{n+2}$ with $\HM{n+1}$). The limitation $n > d/2+1$
was motivated by the basic features of the bilinear map $\PPP$:
this sends continuously $\HM{n} \times \HM{n+1}$ into $\HM{n}$
whenever $n > d/2$ and fulfills the known Kato inequality,
essential for our purposes, if $n > d/2 +1$ (see Section
\ref{notations}; the Kato inequality reviewed therein
reads $| \la \PPP(v, w) | w \ra_n | \leqs G_n \| v \|_n \| w \|^2_n$,
with $G_n$ a suitable constant, also depending on $d$).
\par The method proposed in \cite{appeul} was inspired by
Chernyshenko \emph{et al.} \cite{Che} (and, partly, by \cite{due}
\cite{accau}); given an approximate solution $\ua \in C([0,\Ta),
\HM{n+2}) \cap C^1([0,\Ta), \HM{n})$ of the NS Cauchy problem, it
allows to infer a lower bound on the interval of existence of the
exact NS solution $u$, and an upper bound on the $\HM{n}$ distance
between $u(t)$ and $\ua(t)$. This is obtained via an \emph{a
posteriori} analysis of $\ua$ relying on the norms \beq \| \big({d
\ua \over d t} - \nu \Delta \ua - \PPP(\ua, \ua) - f\big)(t)
\|_n~, \qquad \| \ua(0) - u_0 \|_n \label{norme0} \feq \beq \|
\ua(t) \|_n~, \qquad \| \ua(t) \|_{n+1} \label{norme} \feq ($t \in
[0,\Ta)$), which measure the ``differential error'', the ``datum
error'' and the ``growth'' of $\ua$. The above norms, or some
upper bounds for them, determine some inequalities for an unknown
function $\Rr_n \in C([0,\Tc), \reali)$, that we have called the
\emph{control inequalities}; these consist of a differential
inequality for $\Rr_n$, supplemented with an inequality for the
initial value $\Rr_n(0)$. The main result of \cite{appeul} is the
following: if the control inequalities have a solution $\Rr_n$
with domain $[0,\Tc)$, then the exact solution $u \in C([0,T),
\HM{n+2}) \cap C^1([0,T), \HM{n})$ of the NS Cauchy problem is
such that \beq T \geqs \Tc~, \qquad \| u(t) - \ua(t) \|_n \leqs
\Rr_n(t)~ \mbox{on $[0,\Tc)$}. \label{15} \feq On the other hand,
it is known that the NS Cauchy problem with $C^\infty$ initial
data and forcing has a $C^\infty$ solution. Thus, it can be of
interest to propose a variant of the approach of \cite{appeul}
where the a posteriori analysis of approximate solutions and its
implications on the exact solution are presented in a $C^\infty$
functional setting; this is the aim of the present paper. The
starting points of our analysis are the following ones:
\begin{itemize}
\item[(a)] one can introduce the \Fre space $\HM{\infty}$,
intersection of the finite order Sobolev spaces $\HM{p}$
as $p$ ranges in $\reali$ (or in any subset of $\reali$ unbounded from
above, e.g., $\naturali$). This coincides (algebraically and topologically)
with the space of the $C^\infty$ vector
fields $v: \Td \vain \reali^d$ having zero divergence and
mean.
\item[(b)] The NS bilinear map $\PPP$ fulfills known inequalities
where a norm or an inner product of arbitrarily large Sobolev order $p$
has a bound involving the Sobolev norms of order
$p$ (or $p+1$) and of a fixed, lower order $n$ (or $n+1$):
see the forthcoming Eqs.\,\rref{basineqa}\rref{katineqa}.
These are ``tame'' inequalities in the general sense
used in studies on the Nash-Moser implicit function theorem \cite{Ham}.
\end{itemize}
Under the assumption of an initial datum $u_0 \in \HM{\infty}$ and of a forcing
 $f \in$ $C^\infty([0,+\infty), \HM{\infty})$, and given
an approximate solution
$\ua \in C^1([0,\Ta), \HM{\infty})$
of the NS Cauchy problem, the main results of the paper are as follows:
\begin{itemize}
\item[(i)]
we can start from the norms \rref{norme0}\rref{norme} for a given Sobolev
order $n > d/2+1$ and reconsider
the control inequalities of \cite{appeul} in
an unknown function  $\Rr_n \in C([0,\Tc), \reali)$; if
these have a solution $\Rr_n$ of domain $[0,\Tc)$,
then the exact solution
$u \in C^\infty([0,T), \HM{\infty})$ of
the NS Cauchy problem fulfills the bounds \rref{15}.
If $\Rr_n$ is global ($\Tc = +\infty$), then $u$ is global
as well ($T=+\infty$).
\item[(ii)]
For any $p > n$, the Sobolev norms
\beq {~} \| \big({d \ua \over d t} - \nu \Delta \ua - \PPP(\ua, \ua) - f\big)(t) \|_p~, \qquad
\| \ua(0) - u_0 \|_p~, \label{normep0} \feq
\beq \| \ua(t) \|_p\,,~~\| \ua(t) \|_{p+1} \label{normep} \feq
and the function $\Rr_n$ of item (i) can be used
to construct linear inequalities for an unknown
real function on $[0,\Tc)$; these have an explicitly
computable solution $\Rr_p : [0,\Tc) \vain \reali$, and we prove that
\beq \| u(t) - \ua(t) \|_p \leqs \Rr_p(t) \feq
for all $t \in [0,\Tc)$.
\end{itemize}
In a few words: a suitable a posteriori analysis of $\ua$
gives bounds on the exact NS solution $u$ in
the Sobolev norms of arbitrary order. \par
The simplest application of the above scheme
is set up choosing $\ua(t) := 0$ for all $t \geqs 0$. With
this choice (and assuming $f=0$ as a further simplification),
we can obtain from the control inequalities simple
and fully explicit bounds on the NS solution $u$ and
its time of existence $T$; these imply, for
example, that $u$ is global and exponentially decaying
in all Sobolev norms if the datum is sufficiently small,
to be precise if $\| u_0 \|_n \leqs \nu/G_n$ (a result
which is not at all surprising but might have some
interest in its present, fully quantitative
formulation). As a matter of fact,
the general scheme (i)(ii) outlined before
has been mainly devised for more sophisticated
applications; however, these results will just be sketched in the present paper.
\par
The paper is organized as follows. Section \ref{notations}
and the related Appendix \ref{appeco}
are devoted to some preliminaries; in particular, they describe
the inequalities for the NS bilinear map $\PPP$ which
have been mentioned in the previous discussion. Section \ref{causm} and the related
Appendix \ref{appecau} deal with the NS Cauchy problem in an $\HM{\infty}$ framework.
Section \ref{secmain} contains
the main result of the paper (Proposition \ref{main}), corresponding
to the previous items (i)(ii). Section \ref{appzero}
applies this results with the simple
choice $\ua(t) :=0$ for all $t$ (assuming
$f=0$ for simplicity). Section \ref{appli} indicates
the possibility of more sophisticated choices
of $\ua$, reconsidering from the viewpoint of the present work
some applications of the general
method of \cite{appeul} presented therein and
in some related works \cite{apprey} \cite{bnw} \cite{padova}  \cite{Forth}.
In these applications $\ua$ was a Galerkin
approximant, or a truncated
expansion in the Reynolds number or in time,
typically in dimension $d=3$; here
we only indicate how these applications
could be refined along the scheme of the present paper, leaving
the details to future works.
\vfill \eject \noindent
{~}
\vskip -2cm \noindent
\section{Preliminaries}
\label{notations}
\textbf{Function spaces of vector fields on the torus.}
Throughout this paper we work
on a torus $\Td := (\reali/2 \pi \interi)^d$
of any dimension $d \in \{2,3,...\}$, keeping all the notations employed
in \cite{appeul} (and in most of the other works of ours, cited in the bibliography).
In particular, we write $\mathbb{D}'(\Td) \equiv
\mathbb{D}'$ for the space of $\reali^d$-valued distributions on $\Td$;
each $v \in \mathbb{D}'$ has a weakly convergent
Fourier expansion $v = (2 \pi)^{-d/2} \sum_{k \in \Zd} v_k e^{i k \sc x}$,
with coefficients $v_k = \overline{v_{-k}} \in \complessi^ d$.
The mean value $\la v \ra$ is, by definition, the action of $v$
on the constant test function $(2 \pi)^{-d}$, and $\la v \ra = (2 \pi)^{-d/2} v_0$.
The Laplacian of $v \in \mathbb{D}'$ has Fourier coefficients $(\Delta v)_k := -|k|^2 v_k$;
if $\la v \ra = 0$ and $n \in \reali$, we define $(-\Delta)^{n/2} v$
to be the element of $\mathbb{D}'$ with mean zero and Fourier coefficients $((-\Delta)^{n/2} v)_k =
|k|^n v_k$ for $k \in \Zd \setminus \{ 0 \}$.
\par
Let us consider the space $L^2(\Td, \reali^ d) \equiv \bb{L}^2$,
with its standard inner product $\la~ |~ \ra_{L^2}$. For any $n \in \reali$, we consider the
Sobolev space
\parn
\vbox{
\beq {~} \hspace{-0.5cm} \HM{n}(\Td) \equiv \HM{n} :=
\{ v \in \mathbb{D}'~|~~\dive v = 0,~ \la v \ra = 0, ~
(-\Delta)^{n/2} v \in \bb{L}^2~\}  \label{defhn} \feq
$$ = \{ v \in \mathbb{D}'~|~~k \sc v_k = 0 ~\forall k \in \Zd,~~ v_0 = 0, ~
\sum_{k \in \Zd \setminus \{0\}} |k|^{2 n} | v_k |^2 < + \infty \} $$
}
\noindent
(the lowercase symbols $\ss$, $\zz$ are used to recall the vanishing
of the divergence and of the mean). The above
Sobolev space is equipped with
the inner product and with the induced norm
\beq \la v | w \ra_n := \la (-\Delta)^{n/2} v |  (-\Delta)^{n/2} w \ra_{L^2}
= \! \! \! \! \sum_{k \in \Zd \setminus \{0 \}} |k|^{2 n} \overline{v_k} \sc w_k\, , ~~
\| v \|_n := \sqrt{\la v | v \ra_n}~. \label{definner} \feq
One has $\HM{p} \hookrightarrow \HM{n}$ if $p \geqs n$, where
$\hookrightarrow$ indicates a continuous imbedding (more quantitatively:
$\|~\|_p \geqs \|~\|_n$ if $p \geqs n$).
The vector space
\beq {~} \hspace{-0.5cm} \HM{\infty} := \cap_{p \in \reali} \HM{p} \feq
can be equipped with the topology induced by the family
of all Sobolev norms $\|~\|_p$ ($p \in \reali)$. This space and
its topology do not change if
$\reali$ is replaced with any subset of the reals unbounded from above, e.g., $\naturali$;
the countability of the family of norms $\|~\|_p$ ($p \in \naturali$)
ensures that we have a \Fre topology. \par
For $k \in \naturali \cup \{ \infty \}$ we consider the space
\beq \CM{k}(\Td) \equiv \CM{k} := \{ v \in C^k(\Td,\reali^d)~|~~
\dive v = 0,~ \la v \ra = 0~\}~, \feq
which is a Banach space for $k < \infty$ and a \Fre space for $k=\infty$, when
equipped with the sup norms for all derivatives up to order $k$.
Let $h, k \in \naturali$, $p \in \reali$;
then $\CM{h} \hookrightarrow \HM{p}$ if $h \geqs p$ and, by the Sobolev lemma,
$\HM{p} \hookrightarrow \CM{k}$ if $p > k + d/2$. From these facts one easily infers
\beq \HM{\infty} = \CM{\infty} \feq
(which indicates the equality of the above vector spaces and of
their \Fre topologies).
\vfill \eject \noindent
\textbf{The NS bilinear map, and some inequalities for it.}
We have already introduced the notation
$\PPP$ for the fundamental bilinear map in the NS equations on $\Td$,
see Eq.\,\rref{defp}.
\par
Let $n, p$ be real numbers fulfilling the inequalities written hereafter; it is known that
\beq n > d/2,~ v \in \HM{n},~ w \in \HM{n+1} \quad \Rightarrow \quad \PPP(v,w) \in \HM{n} \label{known} \feq
and that there are constants $K_n$, $G_n$, $K_{p n}$, $G_{p n}$ $\in (0,+\infty)$ such that the following holds:
\beq \| \PPP(v, w) \|_n \leqs K_n \| v \|_n \| w \|_{n+1} \qquad \mbox{for $n > d/2$,
$v \in \HM{n}$, $w \in \HM{n+1}$}~, \label{basineq} \feq
\beq | \la \PPP(v, w) | w \ra_n | \leqs G_n \| v \|_n \| w \|^2_n
\qquad \mbox{for $n > d/2 + 1$, $v \in \HM{n}$,
$w \in \HM{n+1}$}~, \label{katineq} \feq
\beq \| \PPP(v, w) \|_p \leqs {1 \over 2} K_{p n} ( \| v \|_p \| w \|_{n+1} + \| v \|_n \| w \|_{p+1})
\label{basineqa} \feq
$$ \qquad \mbox{for $p \geqs n > d/2$, $v \in \HM{p}$, $w \in \HM{p+1}$}~, $$
\beq | \la \PPP(v, w) | w \ra_p | \leqs
{1 \over 2} G_{p n} (\| v \|_p \| w \|_n + \| v \|_n \| w \|_p)\| w \|_p \label{katineqa} \feq
$$ \mbox{for $p \geqs n > d/2 + 1$, $v \in \HM{p}$, $w \in \HM{p+1}$}~. $$
Note that \rref{basineqa} with $p=n$ implies \rref{basineq}, with $K_n := K_{n n}$;
similarly, \rref{katineqa} with $p=n$ gives \rref{katineq} with $G_n := G_{n n}$.
Statements \rref{known} \rref{basineq} indicate that $\PPP$ maps continuously
$\HM{n} \times \HM{n+1}$ to $\HM{n}$. The fact that these statements
hold for all $n > d/2$ also ensures that $\PPP$ maps continuously $\HM{\infty}
\times \HM{\infty}$ to $\HM{\infty}$. \par
Eq.\,\rref{basineq}
is closely related to the basic norm inequalities about
multiplication in Sobolev space, and \rref{katineq} is due to Kato \cite{Kato};
fully quantitative estimates for the constants $K_n$, $G_n$ therein
have been given in our papers \cite{cog} \cite{cok} where
\rref{basineq} and \rref{katineq} are referred to,
respectively, as the basic and
Kato inequalities for $\PPP$.
Eqs.\,\rref{basineqa} \rref{katineqa} could be
referred to as the generalized basic and Kato inequalities; as mentioned
in the Introduction, they are ``tame'' refinements (in the Nash-Moser sense) of
Eqs.\,\rref{basineq} \rref{katineq}. We remark that inequalities
very similar to \rref{katineqa} are used by Temam in \cite{Tem} and by Beale-Kato-Majda in \cite{BKM};
recently, some analogous inequalities have been proposed by Robinson-Sadowski-Silva
\cite{RSS} as a tool to investigate the putative blow-up of NS solutions.
Explicit expressions for the constants $K_{p n}, G_{p n}$ in \rref{basineqa} \rref{katineqa}
are given in Appendix \ref{appeco} and in \cite{coga}.
\vfill \eject \noindent
\section{The NS Cauchy problem in a smooth framework}
\label{causm}
We are now ready to discuss the NS Cauchy problem in
the framework of the space  $\HM{\infty} = \CM{\infty}$. Let us choose
\beq \nu \in [0,+\infty)\, , \qquad f \in C^\infty([0,+\infty), \HM{\infty}) \, ,
\qquad
 \uz \in \HM{\infty}~. \label{asinf} \feq
\begin{prop}
\textbf{Definition.}
The \textsl{
(incompressible) NS Cauchy problem} with viscosity $\nu$,
forcing $f$ and initial datum $\uz$ is the following:
\beq \mbox{Find}~
u \in C^\infty([0, T), \HM{\infty})~~\mbox{such that}
 ~~{d u \over d t} = \nu \Delta u + \PPP(u,u) + f,~~ u(0) = \uz \label{cau} \feq
(with $T \in (0, + \infty]$, depending on $u$).
If $\nu=0$, this will also be called the ``Euler Cauchy problem'' with
datum $\uz$ and forcing $f$.
\end{prop}
\begin{prop}
\label{procau}
\textbf{Proposition}
With $\nu, f, \uz$ as in \rref{asinf}, the following holds. \parn
(i) Problem \rref{cau} has a unique maximal (i.e., not extendable) solution,
from now on denoted by $u$, with a suitable
domain $[0,T)$. Every solution is a restriction of the maximal one. \parn
(ii) (Beale-Kato-Majda blow up criterion.) Let $u, T$ be as before. If $T < + \infty$,
then $\int_{0}^T d t \| rot u(t) \|_{L^\infty} = + \infty$, whence
$\limsup_{t \vain T^{-}} \| \rot u(t) \|_{L^\infty} = +\infty$. \parn
(iii) The result (ii) implies the following: if $T < + \infty$, then
for each real $n > d/2+1$ one has $\int_{0}^T d t \| u(t) \|_{n} = + \infty$,
whence $\limsup_{t \vain T^{-}} \| u(t) \|_{n} = +\infty$.
\end{prop}
The above proposition is known; it combines results from Kato \cite{Kato}, Temam
\cite{Tem} and Beale-Kato-Majda \cite{BKM} on local existence and blow up
for the Euler equations that, as indicated by the authors themselves, have
simple generalizations to NS equations with arbitrary viscosity;
for more details, we refer to Appendix \ref{appecau}. In the approach
of this Appendix, the main reason for local existence in $\HM{\infty}$ is that local existence
can be established in Sobolev spaces of finite but arbitrarily large order,
on a time interval independent of the order; this idea was first
advanced by Temam \cite{Tem}, on the grounds of a blow up criterion
slightly weaker than the one of Beale-Kato-Majda. \par
For completeness we mention that, in the case $\nu >0$,
statement (ii) can be replaced by a blow up criterion of Giga
\cite{Gig}\cite{Koz2} asserting that $\int_{0}^T d t \| u(t) \|^2_{L^\infty} = + \infty$,
and implying $\int_{0}^T d t \| u(t) \|^2_{n} = + \infty$ for
all $n > d/2$; this is not relevant for our present purposes since
the treatment that we propose for the approximate solutions, related to the
Kato and generalized Kato inequalities \rref{katineq} \rref{katineqa},
relies on the Sobolev norms of orders $> d/2+1$.
\section{Approximate solutions of the NS Cauchy problem and control inequalities}
\label{secmain}
Assuming again $\nu, f, \uz$ as in \rref{asinf}, let us
stipulate what follows.
\begin{prop}
\label{defap}
\textbf{Definition.} \textsl{An \emph{approximate solution} of the problem \rref{cau} is
any map $\ua \in$ $C^1([0, \Ta), \HM{\infty})$,
with $\Ta \in (0,+\infty]$. Given such a function,
we stipulate (i) (ii)}. \par\noindent
(i) The \emph{differential error} of $\ua$ is
\beq e(\ua) := {d \ua \over d t} - \nu \Delta \ua - \PPP(\ua,\ua) - f~
\in C([0,\Ta), \HM{\infty})~;  \label{differr} \feq
the \emph{datum error} is
\beq \ua(0) - \uz \in \HM{\infty}~. \feq
(ii) Let $p \in \reali$. A \emph{differential error estimator},
a \emph{datum error estimator} and a \emph{growth estimator} of order $p$ for $\ua$
are a function $\ep_p \in C([0,\Ta), [0,+\infty))$, a
number $\delta_p \in [0,+\infty)$ and
a function $\Dd_p \in C([0,\Ta), [0,+\infty))$ such that, respectively,
\beq \| e(\ua)(t) \|_p \leqs \ep_p(t)~\mbox{~~for $t \in [0,\Ta)$}~, \label{erest} \feq
\beq \| \ua(0) - \uz \|_p \leqs \delta_p~, \label{uaz} \feq
\beq \| \ua(t) \|_p \leqs \Dd_p(t)~\mbox{~~for $t \in [0,\Ta)$}~.
\label{din} \feq
In particular the function $\ep_p(t) := \| e(\ua)(t) \|_p$,
the number $\delta_p := \| \ua(0) - \uz \|_p$ and
the function $\Dd_p(t) := \| \ua(t) \|_p$ will be called the
\emph{tautological} estimators of order $p$ for the differential error,
the datum error and the growth of $\ua$.
\end{prop}
We note that, according to the previous definition,
an approximate solution $\ua$ is in $\HM{\infty}$
and thus is $C^\infty$ at any instant, but is only required to be $C^1$ in time;
stronger regularity conditions, such as the
assumption that $\ua$ is $C^\infty$ in time, are not
necessary for the sequel. \par
The forthcoming lemma presents an estimate on
the time derivative of the Sobolev distance of any order $p > d/2 + 1$
between the exact solution of the NS Cauchy problem and an approximate solution.
This estimate is the basic tool to establish our main result on approximate solutions,
which is contained in Proposition \ref{main}. Before stating the lemma and
the proposition we introduce the following notations and assumptions, to be used
throughout the section: \parn
\begin{itemize}
\item[(I)] $u \in C^\infty([0,T), \HM{\infty})$ is the maximal solution of the
NS Cauchy problem \rref{cau};
\item[(II)] $\ua \in C^1([0,\Ta), \HM{\infty})$ is any
approximate solution of \rref{cau}. For each $p > d/2 + 1$,
$\ep_p, \delta_p$ and $\Dd_p$ are estimators of order $p$ for the differential error,
the datum error and the growth of $\ua$;
\item[(III)] $K_n, G_n$ and $G_{p n}$ are constants
fulfilling the inequalities \rref{basineq} \rref{katineq} \rref{katineqa},
for all real $n$ and $p$ with the limitations indicated therein; \parn
\item[(IV)] for each function $\varphi : [0,\tau) \vain \reali$
($\tau \in (0,+\infty]$), we use the \emph{right, upper Dini derivative}
$\dd{{d^{+} \varphi \over d t} : [0, \tau) \vain (-\infty,+\infty]}$,
$\dd{t \mapsto {d^{+} \varphi \over d t}(t) := \limsup_{h \vain 0^{+}} {\varphi(t + h) - \varphi(t) \over h}}$.
\end{itemize}
\begin{prop}
\label{lem}
\textbf{Lemma.} Consider the $C^1$ function
\beq u - \ua : [0, \Tw) \vain \HM{\infty}~, \qquad \Tw := \min(T,\Ta)~. \feq
For any real $p$, introduce the norm
$\| u - \ua \|_p : [0, \Tw) \vain [0,+\infty)$,
$t \mapsto \| u(t) - \ua(t)\|_p$
(a continuous function, possibly non-differentiable where $u(t)=\ua(t)$). If
$n, p \in \reali$ are such that $d/2 + 1 < n \leqs p < + \infty$, the following holds
everywhere on $[0, \Tw)$:
\beq {d^{+} \over d t} \, \| u - \ua \|_p \leqs - \nu \| u - \ua \|_p \label{vez1p} \feq
$$  + (G_p \Dd_p + K_p \Dd_{p+1}) \| u - \ua\|_p + G_{p n} \| u - \ua \|_n
\| u - \ua \|_p + \ep_p ~. $$
\end{prop}
\textbf{Proof.} For the sake of brevity, we put
\beq w := u - \ua~. \feq
The function $w$ fulfills
\beq {d w \over d t} = \nu \Delta w + \PPP(\ua, w) + \PPP(w, \ua) + \PPP(w,w) - e(\ua)~,
\label{fo} \feq
\beq w(0) = \uz - \ua(0)~. \label{fp} \feq
Let us consider the function $\| w \|_p$. In a neighborhood of any instant $t$ such that
$w(t) \neq 0$ this function is differentiable, and
\beq {d^{+} \| w \|_p \over d t} =
{d \| w \|_p \over d t} = {1 \over 2 \| w \|_p} {d \| w \|^2_p \over d t} =
{1 \over \| w \|_p} \la {d w \over d t} | w \ra_p \label{dainsp} \feq
$$ = {1 \over \| w \|_p} \big(\nu \la \Delta w | w \ra_p + \la \PPP(\ua, w) | w \ra_p +
\la \PPP(w, \ua) | w \ra_p + \la \PPP(w,w) | w \ra_p - \la e(\ua) | w \ra_p \big)~.
$$
On the other hand, using the Fourier representations for $\Delta$
and for $\la~|~\ra_p$, $\| ~ \|_{p+1}$, $\|~ \|_p$ we easily infer
\beq \la \Delta w | w \ra_p = - \| w \|^2_{p+1} \leqs - \| w \|^2_p~; \label{firp} \feq
moreover, using the inequalities \rref{basineq} \rref{katineq} \rref{katineqa} for $\PPP$,
the Schwarz inequality for $\la~|~\ra_p$ and
the relations \rref{erest} \rref{din} for $\ep_p ,\Dd_p, \Dd_{p+1}$ we get
\beq \la \PPP(\ua, w) | w \ra_p \leqs G_p \| \ua \|_p \| w \|^2_{p} \leqs
G_p \Dd_p \| w \|^2_p~, \feq
\beq \la \PPP(w, \ua) | w \ra_p \leqs \| \PPP(w, \ua) \|_p \| w \|_{p}
\leqs K_p  \| \ua \|_{p+1} \| w \|_{p}^2 \leqs K_p \Dd_{p+1} \|w \|^2_p~, \feq
\beq \la \PPP(w,w) | w \ra_p \leqs G_{p n} \| w \|_n \| w \|^2_p~, \feq
\beq - \la e(\ua) | w \ra_p \leqs \| e(\ua) \|_p \| w \|_p \leqs \ep_p \| w \|_p~.
\label{lasp} \feq
Inserting \rref{firp}-\rref{lasp} into \rref{dainsp} one obtains
the relation
\beq {d^{+} \| w \|_p \over d t} \leqs
- \nu \| w \|_p + (G_p \Dd_p + K_p \Dd_{p+1}) \|w \|_p + G_{p n} \| w \|_n \| w \|_p
+ \ep_p~; \label{ves1p} \feq
this is just the thesis \rref{vez1p}, in a neighborhood of the instant $t$
under consideration for which we were assuming $w(t) \neq 0$. \par
To conclude, we show that Eq.\rref{ves1p}
holds as well at any instant $t$ such that $ w(t) =0$.
In fact, at any such instant we have
\beq {d^{+} \| w \|_p \over d t} \leqs_1 \| {d w \over d t} \|_p =_2
\|  e(\ua) \|_p \leqs \ep_p =_3 \mbox{r.h.s. of \rref{vez1p}}~.
\feq
In the above, the inequality $\leqs_1$ follows from a general
property of the Dini derivative (see, e.g., \cite{Pet}); the
equality $=_{2}$ follows from \rref{fo} and from
$w(t)=0$; the equality $=_3$ uses again the vanishing
of $w$ at this instant. \fine
\begin{rema}
\textbf{Remark.} The inequality \rref{vez1p} was
presented in our work \cite{appeul} in the special case $p=n$
(and in a different framework reviewed in
the Introduction, where $u, \ua$ were
just continuous as maps to $\HM{n+2}$ and $C^1$ as maps to $\HM{n}$);
in this work we pointed out the relations between
this ($p=n$) inequality and a similar result
obtained in \cite{Che}, to which
\cite{appeul} is greatly indebted. The proof
of Lemma \ref{lem} combines ideas from
the cited works with the tame generalization
\rref{katineqa} of the Kato inequality.
\end{rema}
We are now ready to state our main result.
\begin{prop}
\label{main}
\textbf{Proposition.} Consider a real $n > d/2 + 1$, and
assume there is a function $\Rr_n \in C([0,\Tc), \reali)$,
with $\Tc \in (0,\Ta]$, fulfilling the following
\textsl{control inequalities}:
\beq {d^{+} \Rr_n \over d t} \geqs - \nu \Rr_n
+ (G_n \Dd_n + K_n \Dd_{n+1}) \Rr_n + G_n \Rr^2_n + \ep_n
~\mbox{everywhere on $[0,\Tc)$}, \label{cont1} \feq
\beq \Rr_n(0) \geqs \delta_n \label{cont2} \feq
(note that \rref{cont1} \rref{cont2} are fulfilled as equalities
by a unique function in $C^1([0,\Tc),\reali)$, for a suitable $\Tc$).
Then, (i)(ii) hold. \parn
(i) The maximal solution $u$ of the NS Cauchy problem and
its time of existence $T$ are such that
\beq T \geqs \Tc~, \label{tta} \feq
\beq \| u(t) - \ua(t) \|_n \leqs \Rr_n(t) \qquad \mbox{for $t \in [0,\Tc)$}~. \label{furth} \feq
In particular, if $\Rr_n$ is global ($\Tc = +\infty)$ then $u$ is global as well
($T=+\infty$). \parn
(ii) Consider any real $p >n$, and let
$\Rr_p \in C([0,\Tc),\reali)$ be a solution
of the \textsl{linear control inequalities}
\beq {d^{+} \Rr_p \over d t} \geqs - \nu \Rr_p
+ (G_p \Dd_p + K_p \Dd_{p+1} + G_{p n} \Rr_n) \Rr_p + \ep_p
~\mbox{everywhere on $[0,\Tc)$}~, \label{cont1p} \feq
\beq \Rr_p(0) \geqs \delta_p~. \label{cont2p} \feq
Then
\beq \| u(t) - \ua(t) \|_p \leqs \Rr_p(t) \qquad \mbox{for
$t \in [0,\Tc)$}~. \label{urp} \feq
Conditions \rref{cont1p} \rref{cont2p} are fulfilled
as equalities by a unique function $\Rr_p \in$ $C^1([0,\Tc), \reali)$, given explicitly by
\beq \Rr_p(t) = e^{\dd{-\nu t + \Aa_p(t)}}
\Big(\delta_p + \int_{0}^t d s \, e^{\dd{\nu s -\Aa_p(s)}} \ep_p(s) \Big)~\qquad
\mbox{for $t \in [0,\Tc)$}~, \label{rp} \feq
\beq \Aa_p(t) := \int_{0}^ t d s \,
\big(G_p \Dd_p(s) + K_p \Dd_{p+1}(s) + G_{p n} \Rr_n(s)\big)~. \label{ap} \feq
\end{prop}
\fine
\salto
\textbf{Proof.}
(i) We use the inequality \rref{vez1p} of Lemma \ref{lem}
with $p=n$, so that $G_{p n} = G_n$; this inequality reads
\beq {d^{+} \over d t} \, \| u - \ua \|_n \leqs - \nu \| u - \ua \|_n \label{ver1} \feq
$$ + (G_n \Dd_n + K_n \Dd_{n+1}) \| u - \ua\|_n + G_{n} \| u - \ua \|^2_n
+ \ep_n \quad \mbox{on \, $[0, \min(T,\Ta))$}~. $$
Moreover, by the very definition of the estimator $\delta_n$ we have
\beq \| u(0) - \ua(0) \|_n \leqs \delta_n~. \label{ver2} \feq
The inequalities \rref{ver1} \rref{ver2} for $\| u - \ua \|_n$
have the same structure as the control inequalities
\rref{cont1} \rref{cont2} for $\Rr_n$, with the
reverse order relations; now, a standard
comparison theorem \textsl{\`a la} \v{C}aplygin-Lakshmikhantam \cite{Las} \cite{Mitr} gives
\beq \| u(t) - \ua(t) \|_n \leqs \Rr_n(t) \qquad \mbox{for $t \in [0, \min(T, \Ta, \Tc)) =
[0, \min(T, \Tc))$}~. \feq
Finally, one has
\beq \min(T, \Tc) = \Tc~; \feq
in fact, if $T <  \Tc$, for
all $t \in [0,T)$ we would have
$\| u(t) \|_n \leqs \| u(t) - \ua(t) \|_n + \| \ua(t) \|_n \leqs \Rr_n(t) + \Dd_n(t)$
and this would imply $\limsup_{t \to T^{-}} \| u(t) \|_n
\leqs \Rr_n(T) + \Dd_n(T) < +\infty$, contradicting item (iii) of Proposition \ref{procau}.
\parn
(ii) Keeping in mind item (i)
we consider the inequality \rref{vez1p} for
$\| u - \ua \|_p$\,, holding on $[0, \min(T, \Ta))$ and,
\emph{a fortiori}, on the shorter interval $[0, \Tc)$; from
\rref{vez1p} and from $\| u - \ua \|_n \leqs \Rr_n$ we get
\beq {d^{+} \over d t} \, \| u - \ua \|_p \leqs - \nu \| u - \ua \|_p \label{ver1p} \feq
$$  + (G_p \Dd_p + K_p \Dd_{p+1} + \G_{p n} \Rr_n) \| u - \ua\|_p + \ep_p
\quad \mbox{on $[0, \Tc)$}~. $$
We add
to this the relation \rref{uaz} $\| u(0) - \ua(0) \|_p \leqs \delta_p$.
The inequalities \rref{ver1p} \rref{uaz} for $\| u - \ua \|_p$
have the same structure as the inequalities
\rref{cont1p} \rref{cont2p} assumed for $\Rr_p$, with the reverse order relations\,; again,
a comparison theorem \textsl{\`a la} \v{C}aplygin-Lakshmikhantam
gives
\beq \| u(t) - \ua(t) \|_p \leqs \Rr_p(t) \qquad \mbox{for $t \in [0, \Tc)$}~. \feq
Finally, one checks by elementary means that the function $\Rr_p$
defined by \rref{rp}\rref{ap} is the unique $C^1$ function on $[0,\Tc)$ fulfilling conditions \rref{cont1p}\rref{cont2p}
as equalities.
\fine
\begin{rema}
\textbf{Remark.}
Given $\Rr_n$, the linearity of the
control inequalities \rref{cont1p}
for the functions $\Rr_p$ ($p > n$)
is closely related to the linearity
of the inequality \rref{vez1p} with respect
to $\| u - \ua \|_p$; on the other hand,
this feature of \rref{vez1p} depends
essentially on the tame structure
of the generalized Kato inequality \rref{katineqa}.
\end{rema}
\section{A simple application of the previous results}
\label{appzero}
The forthcoming result is an application of Proposition
\ref{main} in which $\ua$ is chosen to be zero at all times
(and the forcing is  assumed to vanish, just for simplicity;
on this point, see the forthcoming Remark \ref{remazero}(ii)).
For $\nu,t \in [0,+\infty)$, let us define
\beq e_{\nu}(t) := \left\{ \barray{ll} \dd{1 - e^{-\nu t} \over \nu} & \mbox{if $\nu > 0$}, \\
t & \mbox{if $\nu = 0$} \farray \right. \label{enu} \feq
(noting that $\lim_{\nu \vain 0^{+}} \dd{1 - e^{-\nu t} \over \nu} = t$). In the
sequel $G_{n}$ and $G_{p n}$ have the usual meaning, see Eqs.\,\rref{katineq}\rref{katineqa};
recall that $G_{p n} = G_n$ for $p=n$.
\begin{prop}
\textbf{Proposition.}
\label{prozero}
Consider the Cauchy problem \rref{cau} for the NS equations
with viscosity $\nu$, zero forcing ($f(t)=0$
for all $t \geqs 0$) and any datum $u_0 \in \HM{\infty}$;
let $u \in C^\infty([0,T), \HM{\infty})$ denote its maximal
solution.
After fixing a real $n > d/2+1$, let us put
\beq \Tc := \left\{ \barray{ll} + \infty & \mbox{if $\nu >0$, $\| u_0 \|_n \leqs {\nu/G_n}$}~, \\
- \dd{1 \over \nu} \log\left(1 - \dd{\nu \over G_n \| u_0 \|_n} \right) & \mbox{if $\nu > 0$,
$\| u_0 \|_n > {\nu/G_n}$,} \\
\dd{1 \over G_n \| u_0 \|_n} & \mbox{if $\nu=0$} \farray \right. \label{ta} \feq
(intending $1/(G_n \| u_0 \|_n) := + \infty$ if $u_0 = 0$).
Then, $u$ and its interval of existence fulfill
\beq T \geqs \Tc~, \label{tcn} \feq
\beq \| u(t) \|_p \leqs
{\| u_0 \|_p \, e^{-\nu t} \over \Big[1 - G_n \| u_0 \|_n e_{\nu}(t)\Big]^{G_{p n}/G_n}}
\qquad \mbox{for all real $p \geqs n$, $t \in [0,\Tc)$}. \label{rnp} \feq
(So, if $\| u_0 \|_n \leqs \dd{\nu / G_n}$,
$u$ is global ($T=+\infty$) and its norm of any Sobolev order decays exponentially
for $t \vain + \infty$.)
\end{prop}
\textbf{Proof.} We consider for the Cauchy problem \rref{cau} the zero approximate
solution
\beq \ua(t) := 0 \qquad \mbox{for all $t \in [0,+\infty)$}~. \feq
The differential error of $\ua$ is zero (since $f=0$), and the datum error
is $\ua(0) - u_0 = -u_0$; so, we have the error and growth estimators
\beq \ep_p(t) := 0\,, \qquad
\delta_p := \| u_0 \|_p \, , \qquad \Dd_p(t) := 0 \qquad \qquad
(p \in \reali)~. \feq
After fixing a Sobolev order $n > d/2+1$, we consider
the control inequalities \rref{cont1} \rref{cont2} corresponding
to these estimators and try to fulfill them as equalities
for an unknown function $\Rr_n \in C^1([0,\Tc), \reali)$;
in this way we are led to the Cauchy problem
\beq {d \Rr_n \over d t} = -\nu \Rr_n + G_n \Rr^2_n~, \qquad \Rr_n(0) = \| u_0 \|_n~. \feq
The maximal solution has domain $[0,\Tc)$ with $\Tc$ as in \rref{ta}, and
is given by
\beq \Rr_n(t) := {\| u_0 \|_n e^{-\nu t} \over 1 - G_n \| u_0 \|_n e_{\nu}(t)} \qquad \mbox{for
$t \in [0,\Tc)$}~. \label{ern} \feq
According to item (i) of Proposition \ref{main} we have
$T \geqs \Tc$ and $\| u(t) \|_n \leqs \Rr_n(t)$ on $[0,\Tc)$; this
justifies Eq.\,\rref{tcn} and gives as well Eq.\,\rref{rnp} for $p=n$
(since the right hand side of this equation
equals $\Rr_n(t)$ when $p=n$).
Now, let $p > n$. Item (ii) of Proposition \ref{main}
with the present estimators gives
\beq \| u(t) \|_p \leqs \Rr_p(t) \qquad \mbox{for
$t \in [0,\Tc)$}~, \label{urpp} \feq
\beq \Rr_p(t) := e^{\dd{-\nu t + \Aa_p(t)}} \| u_0 \|_p~, \qquad
\Aa_p(t) := G_{p n} \int_{0}^ t d s \, \Rr_n(s)~. \label{rpp} \feq
The computation of $\Aa_p$ is elementary, and one concludes
\beq \Rr_p(t) =
{\| u_0 \|_p e^{-\nu t} \over \Big[1 - G_n \| u_0 \|_n e_{\nu}(t)\Big]^{G_{p n}/G_n}}~;
\feq
this result and \rref{urpp} justify Eq.\,\rref{rnp} for $p > n$. \fine
\begin{rema}
\textbf{Remarks}.
\label{remazero}
(i) Proposition \ref{prozero} is an extension
of Proposition 5.2 of \cite{appeul}, where the zero approximate
solution was employed to discuss the NS Cauchy problem
(with zero forcing) in finite order Sobolev spaces; in the cited paper
we obtained Eq.\,\rref{rnp} for the special case $p=n > d/2 + 1$. \parn
\parn
(ii) Proposition \ref{prozero} can be generalized
to the case of nonzero forcing, with suitable assumptions
on it. In this case, using the
zero approximate solution we can again obtain
explicit bounds on the interval of existence
of the exact NS solution $u$ and on its
Sobolev norms; in particular, $u$ is global if $\nu >0$
and $F_n$, $\| u_0 \|_n$ are
sufficiently small for some $n > d/2+1$, where
$F_n := \sup_{t \geqs 0} \| f(t) \|_n$.
\end{rema}
\section{An outline of more sophisticated applications}
\label{appli}
The general framework for approximate NS solutions and
control inequalities devised in \cite{appeul} for
a Sobolev setting of a given finite order has
been employed in the same paper and in the related works \cite{apprey, bnw, padova,Forth}
in a number of applications, typically in dimension
$d=3$, where the following situations have been considered.
\begin{itemize}
\item[(a)] $\nu \geqs 0$ and $\ua$ is a Galerkin approximant;
\item[(b)] $\nu > 0$ and $\ua$ is an expansion in powers of the ``Reynolds number'' $1/\nu$, truncated to
some order $N$. More precisely
$\ua(t) = \sum_{j=0}^{N} (1/\nu^j) u_{(j)}(\nu t)$, where
the coefficients $u_{(j)}$ are determined requiring the differential
error to be $O(1/\nu^N)$, and $\ua(0) = u_0$;
\item[(c)] $\nu=0$ and $\ua$ is an expansion in powers of time, truncated
to some order $N$. More precisely $\ua(t) = \sum_{j=0}^{N} t^j u_j$,
where the coefficients $u_j$ are determined requiring the differential
error to be $O(t^{N})$, and $\ua(0) = u_0$.
\end{itemize}
The strategy (a) or (b) gives a global
solution for the control inequalities of \cite{appeul}
when $\nu$ is above some critical value $\nu_{cr} >0$,
depending on the initial datum; in this situation, one infers that the NS exact
solution $u$ is global as well. \par
The applications presented in the cited works typically
involve initial data in $\HM{\infty}$ (such as the vortices
of Behr-Ne$\check{\mbox{c}}$as-Wu, see \cite{appeul}
\cite{apprey,bnw,padova, Forth},
and the vortices of Taylor-Green
and Kida-Murakami, see \cite{Forth}); the forcing
is often chosen to be zero, and could be assumed in any case
to be in $C^\infty([0,+\infty), \HM{\infty})$. Therefore, such applications
can be reconsidered from the viewpoint of the present
Proposition \ref{main}. \par
The control inequalities \rref{cont1} \rref{cont2}
of a given Sobolev order $n > d/2 + 1$
appearing in this proposition
are in fact identical to the ones of \cite{appeul};
they have been already solved for the
cited applications in that paper and
in \cite{apprey,bnw,padova,Forth}
(typically, for $d=n=3$).
From the viewpoint of the present Proposition \ref{main},
the existence of a solution $\Rr_n \in C([0,\Tc),\reali)$
for these control inequalities (possibly, with $\Tc = +\infty$) ensures that the
NS Cauchy problem has a solution $u \in C^\infty([0,T),\HM{\infty})$
with $T \geqs \Tc$, and that $\| u(t) - \ua(t) \|_n \leqs \Rr_n(t)$
on $[0,\Tc)$. \par
For each one of the cited applications, using the already known
function $\Rr_n$ with item (ii)
of Proposition \ref{main} we could estimate $\| u(t) - \ua(t) \|_p$
for $t \in [0,\Tc)$ and an arbitrarily large
Sobolev order $p$. Presenting here these
implementations of item (ii) would bring
the length of this paper above a reasonable
bound; we plan to return on this subject
in future works. These will refer to the
numerical values of the constants
$G_{p n}$ appearing in item (ii) of Proposition
\ref{main}, obtained on the grounds of Appendix \ref{appeco} and \cite{coga}.
\vskip 0.6cm \noindent
\textbf{Acknowledgments.}
This work was partly supported by INdAM, INFN and by MIUR, PRIN 2010
Research Project  ``Geometric and analytic theory of Hamiltonian systems in finite and infinite dimensions''.
We are grateful to an anonymous referee for some remarks and a bibliographical indication,
that were very useful to improve the presentation of our results.
\vfill \eject \noindent
\appendix
\section{Appendix. On the constants $\boma{K_{p n}, G_{p n}}$ in Eqs.\,\rref{basineqa} \rref{katineqa}}
\label{appeco}
The constants $K_{n}, G_n$ in Eqs.\,\rref{basineq} \rref{katineq} were estimated
in \cite{cog} \cite{cok}. In \cite{coga} the approach of these papers
is extended to the ``tame'' inequalities \rref{basineqa} \rref{katineqa},
and it is shown that the constants therein can be taken as follows:
\beq K_{p n} = {1 \over (2 \pi)^{d/2}} \sqrt{ \sup_{k \in \Zd \setminus \{ 0 \}} \KK_{p n}(k)}~,
\qquad G_{p n} = {1 \over (2 \pi)^{d/2}} \sqrt{ \sup_{k \in \Zd \setminus \{ 0 \}} \GG_{p n}(k)}~, \label{a1} \feq
where $\KK_{p n}, \GG_{p n} : \Zd \setminus \{0 \} \vain (0,+\infty)$ are defined by
\beq \KK_{p n}(k) := 4 |k|^{2 p} \sum_{h \in \Zd \setminus \{ 0, k \}}
{C^2_{h, k} \over (|h|^{p} |k-h|^{n + 1} + |h|^n |k-h|^{p+1})^2}~,  \label{a2} \feq
\beq \GG_{p n}(k) := 4 \sum_{h \in \Zd \setminus \{ 0, k \}}
{(|k|^p - |k - h|^p)^2 \, C^2_{h, k}  \over (|h|^p |k-h|^{n} + |h|^n |k-h|^{p})^2}~. \label{a3} \feq
The coefficient $C_{h, k}$ in the above formulas is any upper bound on
the norm of the bilinear map $h^{\perp} \times (k-h)^{\perp} \vain k^{\perp}$,
$(a, b) \mapsto (k-h) \sc a \, \LP_k b$\,, where
$\perp$ indicates the orthogonal complement in $\complessi^d$ and $\LP_k$
the projection of $\complessi^d$ onto $k^{\perp}$; one can take
\beq C_{h, k} = {|h \wedge k| \over |h|} = |k| \sin \theta(h, k) \label{chk} \feq
where $\wedge$ indicates the exterior product (the usual vector product, if $d=3$)
and $\theta(h, k) \in [0,\pi]$ is the angle between $h$ and $k$. \par
For more details, and for a description of the numerical procedures
to compute $\KK_{p n}$, $\GG_{p n}$ and their sups, we refer to \cite{coga}. For
$p=n$ and $C_{h, k}$ as in \rref{chk}, the expressions \rref{a1} \rref{a2} \rref{a3}
for $K_{p n}$ and $G_{p n}$
agree with the ones proposed in \cite{cog} \cite{cok} for $K_n$ and $G_n$.  Admittedly, we
do not know if the constants determined by \rref{a1} \rref{a2} \rref{a3} are sharp.
\vfill \eject \noindent
\section{Appendix. On the NS Cauchy problem, and
the proof of Proposition \ref{procau}}
\label{appecau}
In this Appendix we stipulate
\beq \nu \in [0,+\infty), \qquad \si : = \left\{ \barray{cc} 1 & \mbox{if $\nu=0$}, \\
2 & \mbox{if $\nu >0$}. \farray \right. \label{nusi} \feq
Our aim is to sketch a proof of Proposition \ref{procau} about the
NS Cauchy problem in a framework based on
$\HM{\infty}$. This relies on known ``hard'' results
on the Cauchy problem in finite order Sobolev
spaces, summarized hereafter.
\begin{prop}
\textbf{Definition.} Let
\beq p \in \reali, \, p > d/2 + 1 \,, \quad f \in C([0,+\infty), \HM{p})\,,
\quad  \uz \in \HM{p}~. \label{asep} \feq
The \textsl{
(incompressible) NS Cauchy problem} with viscosity $\nu$, Sobolev order $p$,
forcing $f$ and initial datum $\uz$ is the following:
\beq \mbox{Find}~
u \in C([0, T), \HM{p}) \cap C^1([0,T), \HM{p - \si}) \quad \mbox{such that} \label{caup} \feq
$$ {d u \over d t} = \nu \Delta u + \PPP(u,u) + f~, \qquad u(0) = \uz $$
(with $T \in (0, + \infty]$, depending on $u$).
If $\nu=0$, this will also be called the ``Euler Cauchy problem''.
\end{prop}
\begin{prop}
\label{procaup}
\textbf{Proposition.} With $\nu, \si, p, f, \uz$ as
in \rref{nusi} \rref{asep}, the following holds. \parn
(i) Problem
\rref{cau} has a unique maximal (i.e., not extendable) solution, from
now on indicated with $u$, with a suitable domain $[0,T)$. Every
solution is a restriction of the maximal one.
\parn
(ii)
(Beale-Kato-Majda blow up criterion.) Let $u, T$ be as before. If
$T < + \infty$, then $\int_{0}^T d s \| \rot u(t) \|_{L^\infty}
= +\infty$, whence
$\limsup_{t \vain T^{-}} \| \rot u(t) \|_{L^\infty} = + \infty$. \parn
(iii) The result (ii) implies the
following: if $T < + \infty$, then for each real $n$ such that $d/2+1 < n \leqs p$ one
has $\int_{0}^T d t \| u(t) \|_{n} = + \infty$, whence $\limsup_{t
\vain T^{-}} \| u(t) \|_{n} = +\infty$.
\end{prop}
\textbf{Proof.} (i) was proved by Kato in \cite{Kat2}; as for (ii), see the
original Beale-Kato-Majda result \cite{BKM} (with its
extension by Kozono and Taniuchi \cite{Koz} to an arbitrary space dimension $d$),
or the book by Majda and Bertozzi
\cite{Mabe}. (iii) follows noting that, by the Sobolev imbedding lemma,
$\| \rot u(t) \|_{L^\infty} \leqs \mbox{constant} \| u(t) \|_n$ for all $t \in
[0,T)$. \fine
For completeness let us also mention that, prior to
\cite{BKM}, Temam \cite{Tem} proved in place of (ii) the
following, slightly weaker result (ii'):  if $T < + \infty$, then $\int_{0}^T d
t\, \| \partial u(t) \|_{L^\infty} = +\infty$, where $\partial u$ stands for
the Jacobian matrix $(\partial_r u_s)$; (ii') and the Sobolev lemma are
sufficient to infer the previous statement (iii). (For $\nu >0$ these blow-up criteria
could be enriched with the Giga criterion \cite{Gig} \cite{Koz2} $\int_{0}^T \| u(t) \|^2_{L^\infty}< + \infty$,
implying $\int_{0}^T d t \| u(t) \|^2_{n} =+\infty$
for $d/2 < n \leqs p$.)
\vfill \eject \noindent
{~}
\vskip -0.8cm
Using Proposition \ref{procaup}, it is easy to
derive Proposition \ref{procau} about NS Cauchy problem in $\HM{\infty}$; the
argument employed hereafter is very similar to one proposed by Temam \cite{Tem}
in the particular case of the Euler equations.
\salto \textbf{Proof of Proposition \ref{procau}.}
Recall the assumption \rref{asinf} $f \in C^\infty([0,+\infty),
\HM{\infty})$, $\uz \in \HM{\infty}$. \parn
\textsl{Step 1. The notations
\rref{caup}$_p$, $u_p, T_p$}. In the sequel we need to consider the Cauchy
problem \rref{caup} for many choices of the Sobolev order $p > d/2+1$. To avoid confusion,
we indicate with \rref{caup}$_p$ this Cauchy
problem and (provisionally) write $u_p$ for its maximal solution, of domain
$[0,T_p)$.
\parn
\textsl{Step 2. Let $p, q > d/2 + 1$,
with $p \leqs q$. Then $T_{q} \leqs T_p$ and
$u_q = u_p \restriction [0, T_q)$.} In fact
$u_q$ is a solution of \rref{caup}$_q$,
which implies that $u_q$ is as well a solution
of \rref{caup}$_p$; on the other hand, any solution of \rref{caup}$_p$
is a restriction of $u_p$.
\parn
\textsl{Step 3. Let $p, q > d/2 + 1$.
Then $T_{q} = T_p$ and $u_q = u_p$.}
It suffices to prove this for $p \leqs q$.
In this case we have the result of Step 2,
and we must just show that
$T_q = T_p$. Assume this does not hold;
then $T_q < T_p$ due to Step 2 and,
in particular, $T_q < +\infty$. So, by item
(iii) of Proposition \ref{procaup}
(with $p, n$ replaced by $q, p$)
we have $\limsup_{t \vain T^{-}_q} \| u_q(t) \|_{p} = +\infty$.
On the other hand, the result of Step 2
implies $\limsup_{t \vain T^{-}_q} \| u_q(t) \|_{p} =
\limsup_{t \vain T^{-}_q} \| u_p(t) \|_{p} = \| u_p(T_q) \|_p
< + \infty$. So we have a contradiction; the conclusion is
$T_q = T_p$.
\parn
\textsl{Step 4. The function $u$.} We now denote with $T$ the common value of
$T_p$ for all $p > d/2+1$, and with $u$ the function $u_p$ for any
such $p$. For any $p > d/2 + 1$, $u$ is continuous from $[0,T)$
to $\HM{p}$; this implies $u \in C([0,T), \HM{\infty})$. \parn
\textsl{Step 5.
$u \in C^\infty([0, T), \HM{\infty})$, and $u$ is a solution of Cauchy problem
\rref{cau}.} For any $p > d/2+1$, the function $u = u_p$ is in
$C^1([0,T), \HM{p-\si})$ and fulfills $d u/d t = \nu \Delta u + \PPP(u,u) + f$
on $[0,T)$, $u(0) = \uz$. By the arbitrariness of $p$, we infer that $u \in
C^1([0,T), \HM{\infty})$ and $u$ fulfills the previous differential
equation in the framework of differential calculus in $\HM{\infty}$. But
$\Delta$ is continuous from $\HM{\infty}$ to itself, $\PPP$ is a continuous
bilinear map from $\HM{\infty} \times \HM{\infty}$ to $\HM{\infty}$ and $f$ is
$C^1$ (in fact $C^\infty$) from $[0,T)$ to $\HM{\infty}$; so, $d u/d t = \nu
\Delta u + \PPP(u,u) + f \in C^1([0,T), \HM{\infty})$ and this implies $u \in
C^2([0,T), \HM{\infty})$. An iteration of this argument shows that $u \in
C^k([0,T), \HM{\infty})$ for each $k \in \naturali$, thus proving that $u \in
C^\infty([0, T), \HM{\infty})$.
\parn
\textsl{Step 6. Let $u' \in C^\infty([0, T'), \HM{\infty})$ be a solution of
the Cauchy problem \rref{cau}; then $T' \leqs T$ and $u' = u \restriction
[0,T')$ (thus, $u$ is the unique not extendable solution of \rref{cau}).} In
fact, for any $p > d/2+1$, $u'$ is as well a solution of the Cauchy
problem \rref{caup}$_p$; thus, by Proposition \ref{procaup} $u'$ is a
restriction of the maximal solution $u_p$ of \rref{caup}$_p$, that coincides
with $u$.
\parn
\textsl{Step 7. Let $T < + \infty$; then $\int_{0}^T d t \| rot u(t) \|_{L^\infty}$
$ = + \infty$, whence $\limsup_{t \vain T^{-}} \| \rot u(t)
\|_{L^\infty}$ $= +\infty$. For each $n > d/2+1$ this implies
$\int_{0}^T d t \| u(t) \|_{n} = + \infty$, whence $\limsup_{t \vain T^{-}}$
$\| u(t) \|_{n} = +\infty$.} The statements on $\| \rot u(t) \|_{L^\infty}$ follow
choosing any $p > d/2+1$ and applying Proposition \ref{procaup} to
the function $u_p=u$. The statements on $\| u(t) \|_n$ follow using again the
Sobolev inequality. \fine

\end{document}